\documentclass[envcountsect,12pt]{llncs}
\usepackage{fancyheadings}
\usepackage{amssymb}
\usepackage{hyperref}
\usepackage[numeric]{amsrefs}
\makeatletter

\renewcommand{\BibLabel}{%
     \hyper@anchorstart{cite.\CurrentBib}\relax\thebib.\hyper@anchorend%
}
\makeatother

\usepackage{amsmath}
\usepackage{amsfonts}
\ifx\pdftexversion\undefined
  \usepackage[dvips]{graphicx}
\else
  \usepackage[pdftex]{graphicx}
\fi

\ifx\pdftexversion\undefined
\else
\fi
\newtheorem{cor}[theorem]{Corollary}
\newtheorem{lem}[theorem]{Lemma}
\newtheorem{prop}[theorem]{Proposition}

\newcommand{\lemref}[1]{Lemma~\ref{#1}}
\newcommand{\secref}[1]{Section~\ref{#1}}

\numberwithin{equation}{section}
\renewcommand\a{\alpha}         \renewcommand\b{\beta}

\renewcommand\d{\delta}
\newcommand\e{\varepsilon}
\renewcommand\l{\lambda}
\renewcommand\L{\Lambda}

\newcommand\G{\Gamma}
\newcommand\f{\frac}
\renewcommand{\div}{\mid}
\newcommand{\ndiv}{\nmid}

\newcommand{\Z}{{\mathbb{Z}}}

\newcommand{\C}{{\mathbb{C}}}

\newcommand{\Q}{{\mathbb{Q}}}

\newcommand{\Cl}{\operatorname{Cl}}
\newcommand{\End}{\operatorname{End}}
\newcommand{\Aut}{\operatorname{Aut}}
\newcommand{\Hom}{\operatorname{Hom}}
\newcommand{\Trace}{\operatorname{Trace}}

\newcommand{\disc}{\operatorname{disc}}

\newcommand{\polylog}{\operatorname{polylog}}
\renewcommand\O{{\mathcal O}}

\newcommand\F{{\mathbb F}}
\renewcommand\H{{\cal H}}
\newcommand{\Graph}{\mathcal{G}}

\renewcommand\i{^{-1}}
\renewcommand\({\left(}
\renewcommand\){\right)}
\newcommand{\jacobi}[2]{{\genfrac{(}{)}{}{}{#1}{#2}}}
\newcommand{\ignore}[1]{{}}
\begin{document}

\title{
Do All Elliptic Curves of the Same Order Have the Same Difficulty of
Discrete Log? }
\author{David Jao\inst{1}, Stephen
Miller\thanks{Partially supported by NSF grant DMS-0301172 and an
Alfred P. Sloan Foundation Fellowship.}\inst{,2,3}, and
Ramarathnam Venkatesan\inst{1}}
\institute{Microsoft Research, 1 Microsoft Way, Redmond WA 98052 \\
\email{\{davidjao,venkie\}@microsoft.com} \and Department of
Mathematics, Rutgers University \\
110 Frelinghuysen Rd, Piscataway, NJ 08854-8019
\\ \and
Einstein Institute of Mathematics, Edmond J. Safra Campus, Givat Ram\\
The Hebrew University of Jerusalem, Jerusalem 91904 Israel \\
\email{miller@math.huji.ac.il} }

\maketitle
\begin{abstract}
The aim of this paper is to justify the common cryptographic
practice of selecting elliptic curves using their order as the
primary criterion.  We can formalize this issue by asking  whether
the discrete log problem (\textsc{dlog}) has the same difficulty
for all curves over a given finite field with the same order.  We
prove that this is essentially true by showing polynomial time
random reducibility of \textsc{dlog} among such curves, assuming
the Generalized Riemann Hypothesis (GRH). We do so by constructing
certain expander graphs, similar to Ramanujan graphs,  with
elliptic curves as nodes and low degree isogenies as edges.
 The result is obtained from the rapid mixing of random walks on this graph. Our
proof works only for curves with (nearly) the same endomorphism
rings. Without this technical restriction such a \textsc{dlog}
equivalence might be false; however, in practice the restriction
may be moot, because all known polynomial time techniques for
constructing equal order curves produce only curves with nearly
equal endomorphism rings.

\vspace{.2cm}

\end{abstract}

\section{Introduction}

\label{introsec}

Public key cryptosystems based on the elliptic curve discrete
logarithm (\textsc{dlog}) problem \cite{miller,koblitz} have received
considerable attention because they are currently the most widely used
systems whose underlying mathematical problem has yet to admit
subexponential attacks (see \cite{blake,washington,handbook}). Hence
it is important to formally understand how the choice of elliptic
curve affects the difficulty of the resulting \textsc{dlog}
problem. This turns out to be more intricate than the corresponding
problem of \textsc{dlog} over finite fields and their selection.

To motivate the questions in this paper, we begin with two
observations.  First, we note that one typically picks an elliptic
curve at random, and examines its group order (e.g.~to check if it is
smooth) to decide whether to keep it, or discard it and pick another
one. It is therefore a natural question whether or not \textsc{dlog}
is of the same difficulty on curves over the same field with the same
number of points. Indeed, it is a theorem of Tate that curves $
E_{1}$ and $E_{2}$ defined over the same finite field $\F_q$ have the
same number of points if and only if they are \emph{isogenous},
i.e.,~there exists a nontrivial algebraic group homomorphism $\phi
\colon E_{1}\rightarrow E_{2}$ between them. If this $\phi$ is
efficiently computable and has a small kernel over $\F_q$, we can
solve \textsc{dlog} on $E_1$, given a \textsc{dlog} oracle for $E_2$.

Secondly, we recall the observation that \textsc{dlog} on
$(\Z/p\Z)^{\ast}$ has \emph{random self-reducibility}: given any
efficient algorithm $A(g^x)=x$ that solves \textsc{dlog} on a
polynomial fraction of inputs, one can solve \emph{any} instance
$y=g^x$ by an expected polynomial number of calls to $A$ with
\emph{random} inputs of the form $A(g^r y)$.  Thus, if \textsc{dlog}
on $(\Z/p\Z)^\ast$ is hard in a sense suitable for cryptography at
all (e.g., has no polynomial on average attack), then all but a
negligible fraction of instances of \textsc{dlog} on $(\Z/p\Z)^\ast$
must necessarily be hard. This result is comforting since for
cryptographic use we need the \textsc{dlog} problem to be hard with
overwhelming probability when we pick inputs at random. The same
random self-reduction statement also holds true for \textsc{dlog} on
any abelian group, and in particular for \textsc{dlog} on a
\emph{fixed} elliptic curve.  We consider instead the following
question: given a polynomial time algorithm to solve \textsc{dlog}
on some positive (or non-negligible) fraction of isogenous elliptic
curves over ${\mathbb{F}}_{q}$, can we solve \textsc{dlog} for {\em
all} curves in the same isogeny class in polynomial time?  In this
paper we show that the answer to this question is essentially yes,
by proving (assuming GRH) the mixing properties of random walks of
isogenies on elliptic curves.  It follows that if \textsc{dlog} is
hard at all in an isogeny class, then \textsc{dlog} is hard for all
but a negligible fraction of elliptic curves in that isogeny class.
This result therefore justifies, in an average case sense, the
cryptographic practice of selecting curves at random within an
isogeny class.

\subsection{Summary of our results}
\label{summary}

The conventional wisdom is that if two elliptic curves over the same
finite field have the same order, then their discrete logarithm
problems are equally hard.  Indeed, this philosophy is embodied in the
way one picks curves in practice.  However, such a widely relied upon
assertion merits formal justification.  Our work shows that this
simplified belief is essentially true for all elliptic curves which
are constructible using present techniques, but with an important
qualification which we shall now describe.

Specifically, let $S_{N,q}$ denote the set of elliptic curves defined
over a given finite field $\F_q$, up to $\bar{\F}_q$-isomorphism, that
have the same order $N$ over $\F_q$. We split $S_{N,q}$ into
\emph{levels} (as in Kohel~\cite{kohel}), where each level represents
all elliptic curves having a particular endomorphism ring over
$\bar{\F}_q$. The curves in each level form the vertices of an {\em
isogeny graph} \cite{mestre,galbraith,volcano}, whose edges represent
prime degree isogenies between curves of degree less than some
specified bound $m$.

\begin{theorem}\label{newmainthm}
(Assuming GRH)  There exists a polynomial $p(x)$, independent
of $N$ and $q$, such that for $m = p(\log q)$ the isogeny graph
$\Graph$ on each level is an expander graph, in the sense that any
random walk on $\Graph$ will reach a subset of size $h$ with
probability at least $\f{h}{2|\Graph|}$ after $\polylog(q)$ steps
(where the implicit polynomial is again independent of $N$ and
$q$).
\end{theorem}

\begin{cor}\label{newmaincor}
(Assuming GRH)  The \textsc{dlog} problem on elliptic curves is
random reducible in the following sense: given any algorithm $A$
that solves \textsc{dlog} on some fixed positive proportion of
curves in fixed level, one can probabilistically solve
\textsc{dlog} on any given curve in that same level with
$\operatorname{polylog}(q)$ expected queries to $A$ with random
inputs.\end{cor}

The proofs are given at the end of \secref{proofs}.  These results
constitute the first formulation of a polynomial time random
reducibility result for the elliptic curve \textsc{dlog} problem
which is general enough to apply to typical curves that one
ordinarily encounters in practice.  An essential tool in our proof
is the \emph{nearly Ramanujan} property of \secref{expander},
which we use to prove the expansion properties of our isogeny
graphs.  The expansion property in turn allows us to prove the
rapid mixing of random walks given by compositions of small degree
isogenies within a fixed level.  Our method uses GRH to prove
eigenvalue separation for these graphs, and provides a new
technique for constructing expander graphs.

The results stated above concern a fixed level.  One might therefore
object that our work does not adequately address the issue of
\textsc{dlog} reduction in the case where two isogenous elliptic
curves belong to different levels. If an attack is {\em balanced},
i.e., successful on each level on a polynomial fraction of curves,
then our results apply. However, if only unbalanced attacks exist,
then a more general equivalence may be false for more fundamental
reasons. Nevertheless, at present this omission is not of much practical
importance. First of all, most random curves over $\F_q$ belong to
sets $S_{N,q}$ consisting of only one level (see \secref{distr}); for
example, in Figure~\ref{standardscurves}, we find that $10$ out of
the $11$ randomly generated curves appearing in international
standards documents have only one level. Second, if the endomorphism
rings corresponding to two levels have conductors whose prime
factorizations differ by quantities which are polynomially smooth,
then one can use the algorithms of \cite{galbraith,kohel} to
navigate to a common level in polynomial time, and then apply
Corollary~\ref{newmaincor} within that level to conclude that
\textsc{dlog} is polynomial-time random reducible between the two
levels.  This situation always arises in practice, because no
polynomial time algorithm is known which even produces a pair of
curves lying on levels whose conductor difference is not
polynomially smooth. It is an open problem if such an algorithm
exists.

\begin{figure}[t]
\begin{center}
\hspace{-.26 cm}
  {\tiny\begin{tabular}[c]{|c|c|c|} \hline
 \textbf{Curve}  &  $\mathbf c_\pi$ \textbf{~(maximal conductor gap in isogeny class)} &
 $\mathbf{P(c_\pi)}$  \textbf{\ = \  largest prime factor of $\mathbf c_\pi$} \\
 \hline NIST P-192 &  $1$ & $1$ \\
\hline
 NIST P-256 &  $3$ & $3$ \\
\hline
 NIST P-384 &  $1$ & $1$ \\
\hline
 NIST P-521 &  $1$ & $1$ \\
\hline
 NIST K-163 &  ${45641 \! \cdot \! 82153\! \cdot \!  56498081 \! \cdot \!  P(c_\pi)}$ & $86110311$ \\
\hline
 NIST K-233 &  $5610641 \! \cdot \! 85310626991\! \cdot \!  P(c_\pi)$ & $150532234816721999$ \\
\hline
 NIST K-283 &  $1697\! \cdot \! 162254089 \! \cdot \!  P(c_\pi)$ & ${ 1779143207551652584836995286271} $ \\
\hline
 NIST K-409 &  ${ 21262439877311 \! \cdot \!  22431439539154506863}\! \cdot \!  P(c_\pi)$ & ${ 57030553306655053533734286593}$ \\
\hline
 &   & ${ 9021184135396238924389891(contd)} $ \\
 NIST K-571  & $3952463 \! \cdot \!  P(c_\pi)$  & ${ 9451926768145189936450898(contd)}$ \\
  &   & ${ 07769277009849103733654828039}$ \\
\hline
 NIST B-163 &  $1$ & $1$ \\
 \hline
 NIST B-233 &  $1$ & $1$ \\
\hline
 NIST B-283 &  $1$ & $1$ \\
\hline
 NIST B-409 &  $1$ & $1$ \\
\hline
 NIST B-571 &  $1$ & $1$ \\
\hline
 IPSec $3^{rd}$OG,$F_{2^{155}}$ &  $1$ & $1$ \\
\hline
 IPSec 4$^{th}OG,$ $F_{2^{185}}$ &  $1$ & $1$ \\
\hline  \end{tabular} }
\caption{A table of curves recommended as international standards
\cite{standardsdocumentsnist,standardsdocumentsipsec}.  Note that the
value of $c_\pi$ for each of the standards curves is small (at most
3), except for the curves in the NIST K (Koblitz curve) family.  These
phenomena are to be expected and are explained in \secref{distr}. Any
curve with $c_\pi = 1$ has the property that its isogeny class
consists of only one level. It follows from the results of
\secref{summary} that randomly generated elliptic
curves with $c_\pi = 1$ (or, more generally, with smooth $c_\pi$) will
have discrete logarithm problems of typical difficulty amongst all
elliptic curves in their isogeny class.
 \label{standardscurves} }
\end{center}
\end{figure}

Our use of random walks to reach large subsets of the isogeny
graph is crucial, since constructing an isogeny between two
specific curves\footnote{If one uses polynomial size circuits
(i.e., polynomial time algorithms with exponential time
pre-processing) for reductions, then one can relate \textsc{dlog}
on two given curves.  This claim follows using the smallness of
diameter of our graphs and the smoothness of the degrees of
isogenies involved. We omit the details.} is believed to be
inherently hard, whereas constructing an isogeny from a fixed
curve to a subset constituting a positive (or polynomial) fraction
of the isogeny graph is proved in this paper to be easy. Kohel
\cite{kohel} and Galbraith \cite{galbraith} present exponential
time algorithms (and thus exponential time reductions) for
navigating between two nodes in the isogeny graph, some of which
are based on random walk \emph{heuristics} which we prove here
rigorously. Subsequent papers on Weil descent
attacks~\cite{ghs,mtw} and elliptic curve trapdoor
systems~\cite{teske} also use isogeny random walks in order to
extend the GHS Weil descent attack~\cite{gauhs} to elliptic curves
which are not themselves directly vulnerable to the GHS attack.
Our work does not imply any changes to the deductions of these
papers, since they also rely on the above heuristic assumptions
involving exponentially long random walks. In our case, we achieve
polynomial time instead of exponential time reductions; this is
possible since we keep one curve  fixed, and  random reducibility
requires only that the other curve be randomly distributed.

\section{Preliminaries}
\label{prelim}

Let $E_1$ and $E_2$ be elliptic curves defined over a finite field
$\F_q$ of characteristic $p$. An isogeny $\phi\colon E_1 \to E_2$
defined over $\F_q$ is a non-constant rational
map defined over $\F_q$ which is also a group homomorphism from
$E_1(\F_q)$ to~$E_2(\F_q)$~\cite[\S III.4]{silverman}. The degree of
an isogeny is its degree as a rational map.
For any elliptic curve $E\colon y^2 + a_1 xy + a_3 y = x^3 + a_2 x^2 +
a_4 x + a_6$ defined over $\F_q$, the Frobenius endomorphism is the
isogeny $\pi\colon E \to E$ of degree $q$ given by the equation
$\pi(x,y) = (x^q,y^q)$. It satisfies the equation
$$
\pi^2 - \Trace(E) \pi + q = 0,
$$ where $\Trace(E) = q + 1 - \#E(\F_q)$ is the trace of the Frobenius
endomorphism of $E$ over $\F_q$. The polynomial $p(X) := X^2 -
\Trace(E) X + q$ is called the characteristic polynomial of $E$.

An endomorphism of $E$ is an isogeny $E \to E$
defined over the algebraic closure $\bar{\F}_q$ of $\F_q$. The
set of endomorphisms of $E$ together
with the zero map forms a ring under the operations of pointwise
addition and composition; this ring is called the endomorphism ring of
$E$ and denoted $\End(E)$. The ring $\End(E)$ is isomorphic either to
an order in a quaternion algebra or to an order in an imaginary
quadratic field~\cite[V.3.1]{silverman}; in the first case we say $E$
is supersingular and in the second case we say $E$ is ordinary.  In
the latter situation, the Frobenius endomorphism $\pi$ can be regarded
as an algebraic integer which is a root of the characteristic
polynomial.

Two elliptic curves $E_1$ and $E_2$ defined over $\F_q$ are said to be
isogenous over $\F_q$ if there exists an isogeny $\phi\colon E_1 \to
E_2$ defined over $\F_q$.  A theorem of Tate states that two curves
$E_1$ and $E_2$ are isogenous over $\F_q$ if and only if $\#E_1(\F_q)
= \#E_2(\F_q)$~\cite[{\S}3]{tate1}.  Since every isogeny has a dual
isogeny~\cite[III.6.1]{silverman}, the property of being isogenous
over $\F_q$ is an equivalence relation on the finite set of
$\bar{\F}_q$-isomorphism classes of elliptic curves defined over
$\F_q$.  We define an isogeny class to be an equivalence class of
elliptic curves, up to $\bar{\F}_q$-isomorphism, under this
equivalence relation; the set $S_{N,q}$ of~\secref{summary} is thus
equal to the isogeny class of elliptic curves over $\F_q$ having
cardinality $N$.

Curves in the same isogeny class are either all supersingular or
all ordinary.  We assume for the remainder of this paper that we
are in the \emph{\bf{ordinary case}}, which is the more
interesting case from the point of view of cryptography in light
of the MOV attack \cite{mov}. Theorem~\ref{newmainthm} in the
supersingular case was essentially known earlier by results of
Pizer \cite{pi90,piz98}, and a proof has been included for
completeness in Appendix~\ref{ss-sec}.

The following theorem describes the structure of elliptic curves
within an isogeny class from the point of view of their
endomorphism rings.

\begin{theorem}\label{isogeny-orders}
Let $E$ and $E'$ be ordinary elliptic curves defined over $\F_q$
which are isogenous over $\F_q$. Let $K$ denote the imaginary
quadratic field containing $\End(E)$, and write $\O_K$ for the
maximal order (i.e.,~ring of integers) of $K$.
\begin{enumerate}
\item The order $\End(E)$ satisfies the property
$\Z[\pi]
  \subseteq \End(E) \subseteq \O_K$. \label{maxorder-containment}
\vspace{.1cm} \item The order
  $\End(E')$ also satisfies $\End(E') \subset K$ and $\Z[\pi]
  \subseteq \End(E') \subseteq \O_K$. \label{embedding-containment}
\vspace{.1cm} \item
The following are equivalent:
\label{isogeny-tfae}
\begin{enumerate}
\vspace{.1cm} \item $\End(E) = \End(E')$. \label{isogeny-tfae-1}
 \vspace{.1cm} \item
 There
exist two isogenies $\phi\colon E \to E'$ and $\psi\colon E \to
E'$ of relatively prime degree, both defined over $\F_q$.
\label{isogeny-tfae-2} \vspace{.1cm} \item $[\O_K:\End(E)] =
[\O_K:\End(E')]$. \label{isogeny-tfae-3} \vspace{.1cm}  \item
$[\End(E):\Z[\pi]] = [\End(E'):\Z[\pi]]$. \label{isogeny-tfae-4}
\end{enumerate}
\vspace{.1cm} \item Let $\phi\colon E \to E'$ be an isogeny from
$E$ to $E'$ of prime degree $\ell$, defined over $\F_q$. Then
either $\End(E)$ contains $\End(E')$ or $\End(E')$ contains
$\End(E)$, and the index of the smaller in the larger divides
$\ell$. \label{order-containment} \vspace{.1cm} \item Suppose
$\ell$ is a prime that divides one of $[\O_K:\End(E)]$ and
$[\O_K:\End(E')]$, but not the other. Then every isogeny
$\phi\colon E \to E'$ defined over $\F_q$ has degree equal to a
multiple of $\ell$. \label{isogeny-degree-constraint}
\end{enumerate}
\end{theorem}
\begin{proof}
\cite[{\S}4.2]{kohel}.
\end{proof}

For any order $\O \subseteq \O_K$, the conductor of $\O$ is
defined to be the integer $[\O_K:\O]$. The field $K$ is called the
CM field of $E$. We write $c_E$ for the conductor of $\End(E)$ and
$c_\pi$ for the conductor of $\Z[\pi]$. Note that this is not the
same thing as the  arithmetic conductor of an elliptic
curve~\cite[\S C.16]{silverman}, nor is it related to the
conductance of an expander graph \cite{sinclairg}. It follows
from~\cite[(7.2) and (7.3)]{cox} that $\End(E) = \Z + c_E \O_K$
and $D = c_E^2 d_K,$ where $D$ (respectively, $d_K$) is the
discriminant of the order $\End(E)$ (respectively, $\O_K$).
Furthermore, the characteristic polynomial $p(X)$ has discriminant
$d_\pi = \disc(p(X)) = \Trace(E)^2 - 4q = \disc(\Z[\pi]) = c_\pi^2
d_K$, with $c_\pi = c_E \cdot [\End(E):\Z[\pi]]$.

Following~\cite{volcano} and~\cite{galbraith}, we say that an isogeny
$\phi\colon E \to E'$ of prime degree $\ell$ defined over $\F_q$ is
``down'' if $[\End(E):\End(E')] = \ell$, ``up'' if
$[\End(E'):\End(E)] = \ell$, and ``horizontal'' if $\End(E) =
\End(E)$. The following theorem classifies the number of degree
$\ell$ isogenies of each type in terms of the Legendre symbol
$\jacobi{D}{\ell}$.

\begin{theorem}\label{direction-classification}
Let $E$ be an ordinary elliptic curve over $\F_q$, with
endomorphism ring $\End(E)$ of discriminant $D$. Let $\ell$ be a
prime different from the characteristic of $\F_q$.
\begin{itemize}
\item Assume $\ell \ndiv c_E$. Then there are
exactly $1 +
  \jacobi{D}{\ell}$ horizontal isogenies $\phi\colon E \to E'$ of degree
  $\ell$.
\begin{itemize}
\vspace{.1cm} \item If $\ell \ndiv c_\pi$, there are no other
isogenies $E \to E'$ of degree $\ell$ over $\F_q$.
\vspace{.1cm} \item If $\ell \div c_\pi$, there are $\ell -
\jacobi{D}{\ell}$ down isogenies of degree $\ell$.
\end{itemize}
\vspace{.1cm} \item Assume $\ell \div c_E$. Then there is one up
isogeny $E \to E'$ of degree $\ell$.
\begin{itemize}
\vspace{.1cm} \item If $\ell \ndiv \f{c_\pi}{c_E}$, there are no
other isogenies $E \to E'$ of degree $\ell$ over $\F_q$.
\vspace{.1cm} \item If $\ell \div \f{c_\pi}{c_E}$, there are
$\ell$ down isogenies of
  degree $\ell$.
\end{itemize}
\end{itemize}
\end{theorem}
\begin{proof}
\cite[{\S}2.1]{volcano} or \cite[{\S}11.5]{galbraith}.
\end{proof}

It follows that the maximal conductor difference between levels in
an isogeny class is achieved between a curve at the top level
(with $\operatorname{End}(E)=\O_K$) and a curve at the bottom
level (with $\operatorname{End}(E)=\Z[\pi]$).

\subsection{Isogeny Graphs}
\label{graph-specification}

We define two curves $E_1$ and $E_2$ in an isogeny class $S_{N,q}$ to
have the same level if $\End(E_1) = \End(E_2)$.  An \emph{isogeny
graph} is a graph whose nodes consist of all elements in $S_{N,q}$
belonging to a fixed level. Note that a horizontal isogeny always goes
between two curves of the same level; likewise, an up isogeny enlarges
the size of the endomorphism ring and a down isogeny reduces the size.
Since there are fewer elliptic curves at higher levels than at lower
levels, the collection of isogeny graphs under the level
interpretation visually resembles a ``pyramid'' or a
``volcano''~\cite{volcano}, with up isogenies ascending the structure
and down isogenies descending.

As in~\cite[Prop.~2.3]{gr1}, we define two isogenies $\phi\colon
E_1 \to E_2$ and $\phi'\colon E_1 \to E_2$ to be equivalent if
there exists an automorphism $\alpha \in \Aut(E_2)$ (i.e., an
invertible endomorphism) such that $\phi' = \alpha \phi$. The
edges of the graph consist of equivalence classes of isogenies
over $\F_q$ between elliptic curve representatives of nodes in the
graph, which have prime degree less than the bound $(\log q)^{2+\d}$
for some fixed constant $\d>0$. The degree bound must
be small enough to permit the isogenies to be computed, but large
enough to allow the graph to be connected and to have the rapid
mixing properties that we want. We will show in
Section~\ref{proofs} that there exists a constant $\d>0$ for which
a bound of $(\log q)^{2+\d}$ satisfies all the requirements,
provided that we restrict the isogenies to a single level.

Accordingly, fix a level of the isogeny class, and let $\End(E) = \O$
be the common endomorphism ring of all of the elliptic curves
in this level. Denote by $\Graph$ the regular graph whose vertices
are elements of $S_{N,q}$ with endomorphism ring $\O$, and whose
edges are equivalence classes of horizontal isogenies defined over
$\F_q$ of prime degree $\le (\log q)^{2+\d}$. By standard facts
from the theory of complex multiplication~\cite[{\S}10]{cox}, each
invertible ideal $\mathfrak{a} \subset \O$ produces an elliptic
curve $\C/\mathfrak{a}$ defined over some number field $L \subset \C$
(called the ring class field of $\O$) \cite[{\S}11]{cox}.  The
curve $\C/\mathfrak{a}$ has complex multiplication by $\O$, and
two different ideals yield isomorphic curves if and only if they
belong to the same ideal class. Likewise, each invertible ideal
$\mathfrak{b} \subset \O$ defines an isogeny $\C/\mathfrak{a} \to
\C/\mathfrak{ab}^{-1}$, and the degree of this isogeny is the norm
$N(\mathfrak{b})$ of the ideal $\mathfrak{b}$. Moreover, for any
prime ideal $\mathfrak{P}$ in $L$ lying over $p$, the reductions
mod $\mathfrak{P}$ of the above elliptic curves and isogenies are
defined over $\F_q$, and every elliptic curve and every horizontal
isogeny in $\Graph$ arises in this way
(see~\cite[{\S}3]{galbraith} for the $p > 3$ case, and~\cite{ghs}
for the small characteristic case).  Therefore, the isogeny graph
$\Graph$ is isomorphic to the corresponding graph $\H$ whose nodes
are elliptic curves $\C/\mathfrak{a}$ with complex multiplication
by $\O$, and whose edges are complex analytic isogenies
represented by ideals $\mathfrak{b} \subset \O$ and subject to the
same degree bound as before. This isomorphism preserves the
degrees of isogenies, in the sense that the degree of any isogeny
in $\Graph$ is equal to the norm of its corresponding ideal
$\mathfrak{b}$ in $\H$.

The graph $\H$ has an alternate description as a Cayley graph on
the ideal class group $\Cl(\O)$ of $\O$. Indeed, each node of $\H$
is an ideal class of $\O$, and two ideal classes
$[\mathfrak{a}_1]$ and $[\mathfrak{a}_2]$ are connected by an edge
if and only if there exists a prime ideal $\mathfrak{b}$ of norm
$\le (\log q)^{2+\d}$ such that $[\mathfrak{a}_1 \mathfrak{b}] =
[\mathfrak{a}_2]$. Therefore, the graph $\H$ (and hence the graph
$\Graph$) is isomorphic to the Cayley graph of the group $\Cl(\O)$
with respect to the generators $[\mathfrak{b}] \in \Cl(\O)$, as
$\mathfrak{b}$ ranges over all prime ideals of $\O$ of norm $\le
(\log q)^{2+\d}$.

\begin{remark}
The isogeny graph $\Graph$ consists of objects defined over the finite
field $\F_q$, whereas the objects in the graph $\H$ are defined over
the number field $L$. One passes from $\H$ to $\Graph$ by taking
reductions mod $\mathfrak{P}$, and from $\Graph$ to $\H$ by using
Deuring's Lifting
Theorem~\cite{deuring,galbraith,langellipticfunctions}. There is no
known polynomial time or even subexponential time algorithm for
computing the isomorphism between $\Graph$ and
$\H$~\cite[{\S}3]{galbraith}. For our purposes, such an explicit
algorithm is not necessary, since we only use the complex analytic
theory to prove abstract graph-theoretic properties of $\Graph$.
\end{remark}

\begin{remark}
The isogeny graph $\Graph$ is typically a symmetric graph, since
each isogeny $\phi$ has a unique dual isogeny $\hat{\phi}\colon
E_2 \to E_1$ of the same degree as $\phi$ in the opposite
direction~\cite[\S III.6]{silverman}. (From the viewpoint of $\H$,
an isogeny represented by an ideal $\mathfrak{b} \subset \O$ has
its dual isogeny represented simply by the complex conjugate
$\bar{\mathfrak{b}}$.)  However, the definition of equivalence of
isogenies from~\cite{gr1} given in \ref{graph-specification}
contains a subtle asymmetry which can sometimes render the graph
$\Graph$ asymmetric in the supersingular case
(Appendix~\ref{ss-sec}). Namely, if $\Aut(E_1)$ is not equal to
$\Aut(E_2)$, then two isogenies $E_1 \to E_2$ can sometimes be
equivalent even when their dual isogenies are not. For ordinary
elliptic curves within a common level, the equation $\End(E_1) =
\End(E_2)$ automatically implies $\Aut(E_1) = \Aut(E_2)$, so the
graph $\Graph$ is always symmetric in this case. Hence, we may
regard $\Graph$ as undirected and apply known results about
undirected expander graphs (as in the following section) to
$\Graph$.
\end{remark}

\section{Expander Graphs}
\label{expander}
  Let $G=({\mathcal{V},E})$ be a finite graph on $h$
vertices $\mathcal{V}$ with undirected edges $\mathcal{E}$.
Suppose $G$ is a regular graph of degree $k$, i.e.,~exactly $k$
edges meet at each vertex. Given a labeling of the vertices
${\mathcal{V}}=\{v_1,\ldots,v_h\}$, the adjacency matrix of $G$ is
the symmetric $h\times h$ matrix $A$ whose $ij$-th entry
$A_{ij}=1$ if an edge exists between $v_i$ and $v_j$, and $0$
otherwise.

It is convenient to identify functions on $\mathcal{V}$ with
vectors in ${\mathbb{R}}^h$ via this labeling, and therefore also
think of $A$ as a self-adjoint operator on $L^2({\mathcal{V}})$.
All of the eigenvalues of $A$ satisfy the bound $|\l |\le k$.
Constant vectors are eigenfunctions of $A$ with eigenvalue $k$,
which for obvious reasons is called the trivial eigenvalue $\l
_{\hbox{triv}}$. A family of such graphs $G$ with $ h\rightarrow
\infty$ is said to be a sequence of \emph{expander graphs} if all
other eigenvalues of their adjacency matrices are bounded away
from $\l _{\hbox{triv}}=k$ by a fixed amount.\footnote{Expansion
is usually phrased in terms of the number of neighbors of subsets
of $G$, but the spectral condition here is equivalent for
$k$-regular graphs and also more useful for our purposes.} In
particular, no other eigenvalue is equal to $k$; this implies the
graph is connected. A \emph{Ramanujan graph} \cite{LPS} is a
special type of expander which has $|\l |\le 2 \sqrt{k-1}$ for any
nontrivial eigenvalue which is not equal to $-k$ (this last
possibility happens if and only if the graph is bipartite). The
supersingular isogeny graphs in Appendix~\ref{ss-sec} are
sometimes Ramanujan, while the ordinary isogeny graphs in
\secref{graph-specification}  do not qualify, partly because their
degree is not bounded. Nevertheless, they still share the most
important properties of expanders as far as our applications are
concerned. In particular their degree $k$ grows slowly (as a
polynomial in $\log|{\mathcal{V}}|$), and they share a
qualitatively similar eigenvalue separation: instead the
nontrivial eigenvalues $\l$ can be arranged to be $O(k^{1/2+\e})$
for any desired value of $\e>0$. Since our goal is to establish a
polynomial time reduction, this enlarged degree bound is natural,
and in fact necessary for obtaining expanders from {\em abelian}
Cayley graphs \cite{alonroich}.  Obtaining {\em any} nontrivial
{\em exponent $\beta<1$} satisfying $\l = O(k^\beta)$  is a key
challenge for many applications, and accordingly we shall focus on
a type of graphs we call ``nearly Ramanujan'' graphs: families of
graphs whose nontrivial eigenvalues $\l$ satisfy that bound.

 A fundamental use of
expanders is to prove the rapid mixing of the random walk on
$\mathcal{V}$ along the edges $\mathcal{E}$. The following rapid
mixing result is standard but we present it below for convenience.
For more information, see \cite{lubotzkybook,sarnakbook,vallette}.

\begin{prop}\label{expander-prop}
\label{rapmix} Let $G$ be a regular graph of degree $k$ on $h$
vertices. Suppose that the eigenvalue $\l $ of any nonconstant
eigenvector satisfies the bound $|\l |\le c$ for some $c<k$. Let
$S$ be any subset of the vertices of $G$, and $x$ be any vertex in
$G$. Then a random walk of any length at least $
\frac{\log{2h/|S|^{1/2}}}{\log{k/c}}$ starting from $x$ will land
in $S$ with probability at least
$\frac{|S|}{2h}=\frac{|S|}{2|G|}$.
\end{prop}
\begin{proof}
There are $k^r$ random walks of length $r$ starting from $x$.  One
would expect in a truly random situation that roughly
$\f{|S|}{h}k^r$ of these land in $S$.  The lemma asserts that for
$r\ge \f{\log{2h/|S|^{1/2}}}{\log{k/c}}$ at least half that number
of walks in fact do.  Denoting the characteristic functions of $S$
and $\{x\}$ as $\chi_S$ and $\chi_{\{x\}}$, respectively, we count
that
\begin{equation}\label{numwalks}
   \#\,\{\text{walks of length $r$ starting at $x$ and landing in $S$}\} \ \ = \ \
   \langle \, \chi_S  \, , \,  A^r   \chi_{\{x\}}  \,  \rangle \, ,
\end{equation}
where $\langle\cdot,\cdot\rangle$ denotes the inner product of
functions in $L^2({\mathcal V})$.  We estimate this as follows.
Write the orthogonal decompositions of $\chi_S$ and $\chi_{\{x\}}$
as
\begin{equation}\label{orthdec}
    \chi_S \ \  = \ \ \f{|S|}{h}\,{\mathbf 1} \ + \ u \ \ \  \ \ \text{and} \ \
    \ \  \ \ \chi_{\{x\}} \ \  = \ \ \f{1}{h}\,{\mathbf 1} \ + \ w\ ,
\end{equation}
where ${\mathbf 1}$ is the constant vector and $\langle u,{\mathbf
1}\rangle = \langle w,{\mathbf 1}\rangle  = 0$.  Then
(\ref{numwalks}) equals the expected value of $\f{|S|}{h}k^r$,
plus the additional term $\langle u,A^r w\rangle$, which is
bounded by $\|u\| \,\|A^r w\|$.  Because $w\perp{\mathbf 1}$ and
the symmetric matrix $A^r$ has spectrum bounded by $c^r$ on the
span of such vectors,
\begin{equation}\label{norms}
    \|u\| \,\|A^r w\| \ \ \le \ \  c^r\,\|u\| \,\| w\| \ \ \le \ \
    c^r  \, \|\chi_S  \|\, \|  \chi_{\{x\}} \| \ \ = \ \
    c^r\,|S|^{1/2}\,.
\end{equation}  For our values of $r$ this is at most half of
$\f{|S|}{h}k^r$, so indeed at least $\f 12 \f{|S|}{h}k^r$ of the
paths terminate in $S$ as was required.
\end{proof}

In our application the quantities $k$, $\frac {k}{k-c}$, and
$\frac{h}{|S|}$ will all be bounded by polynomials in $\log(h)$.
Under these hypotheses, the probability  is at least $1/2$ that
some $\operatorname{polylog}(h)$ trials of random walks of
$\operatorname{polylog}(h)$ length starting from $x$ will reach
$S$ at least once. This mixing estimate is the source of our
polynomial time random reducibility (Corollary~\ref{newmaincor}).

\section{Spectral Properties of the Isogeny Graph}
\label{proofs}

\subsection{Navigating the Isogeny Graph}

Let $\Graph$ be as in Section~\ref{graph-specification}. The isogeny
graph $\Graph$ has exponentially many nodes and thus is too large to
be stored.  However, given a curve $E$ and a prime $\ell$, it is
possible to efficiently compute the curves which are connected to $E$
by an isogeny of degree $\ell$.  These curves $E^\prime$ have
$j$-invariants which can be found by solving the modular polynomial
relation $\Phi_{\ell}(j(E),j(E^\prime))=0$; the cost of this step is
$O(\ell^{3})$ field operations~\cite[11.6]{galbraith}. Given the
$j$-invariants, the isogenies themselves can then be obtained using
the algorithms of~\cite{volcano} (or~\cite{lercier,lm} when the
characteristic of the field is small).  In this way, it is possible to
navigate the isogeny graph locally without computing the entire
graph. We shall see that it suffices to have the degree of the
isogenies in the graph be bounded by $(\log q)^{2+\d}$ to assure the
Ramanujan properties required for $\Graph$ to be an expander.

\subsection{$\protect\theta$-Functions and Graph Eigenvalues}
\label{thetsubsec}

The graph $\H$ (and therefore also the  isomorphic graph $\Graph$)
has one node for each ideal class of $\O$. Therefore, the total
number of nodes in the graph $\Graph$ is the ideal class number of
the order $\O$, and the vertices ${\mathcal{V}}$ can be identified
with ideal class representatives $\{\a_1,\ldots,\a_h\}$.  Using
the isomorphism between $\Graph$ and $\H$, we see that the
generating function $\sum M_{\a_i,\a_j}(n) q^n$ for degree $n$
isogenies between the vertices $\a_i$ and $\a_j$ of $\Graph$ is
given by

\begin{equation}  \label{genfunc}
\sum_{n\,=\,1}^\infty \, M_{\a_i,\a_j}(n) \, q^n \ \ := \ \
\frac{1}{e} \sum_{z\, \in\, \a_i\i\a_j}q^{N(z)/N(\a_i\i\a_j)}\,,
\end{equation}
where $e$ is the number of units in $\O$ (which always equals $2$
for $\disc(\O) > 4$). The sum on the righthand side depends only
on the ideal class of the fractional ideal $\a_i\i\a_j$; by
viewing the latter as
 a lattice in $\C$, we see that $N(z)/N(\a_i\i\a_j)$ is a quadratic
form of discriminant $D$ where $D :=
\disc(\O)$~\cite[p.~142]{cox}. That means this sum is a
$\theta$-series, accordingly denoted as $\theta_{\a_i\i\a_j}(q)$.
It is a holomorphic modular form of weight 1 for the congruence
subgroup $\G_0(|D|)$ of $SL(2,\Z)$, transforming according to the
character $\jacobi{D}{\cdot}$ (see \cite[Theorem
10.9]{iwaniec-blue}).

Before discussing exactly which degrees of isogenies to admit into
our isogeny graph $\Graph$, let us first make some remarks about
the simpler graph on ${\mathcal{V}} = \{\a_1,\ldots,\a_h\}$ whose
edges represent isogenies of degree exactly equal to $n$. Its
adjacency matrix is of course the $h\times h$ matrix
$M(n)=\left[M_{\a_i,\a _j}(n)\right]_{\{1 \le i,j\le h\}}$ defined
by series coefficients in (\ref{genfunc}).  It can be naturally
viewed as an operator which acts on functions on ${\mathcal V}=
\{\a_1,\ldots,\a_h\}$, by identifying them with $h$-vectors
according to this labeling.  We will now simultaneously
diagonalize all $M(n)$, or what amounts to the same, diagonalize
the matrix $A_q=\sum_{n\ge 1}M(n) q^n$ for any value of $q<1$
(where the sum converges absolutely). The primary reason this is
possible is that for each fixed $n$ this graph is an abelian
Cayley graph on the ideal class group $\operatorname{Cl}(\O)$,
with generating set equal to those classes $\a_i$ which represent
an $n$-isogeny.  The eigenfunctions of the adjacency matrix of an
abelian Cayley graph are always given by characters of the group
(viewed as functions on the graph), and their respective
eigenvalues are sums of these characters over the generating set.
This can be seen directly in our circumstance as follows.  The
$ij$-th entry of $A_q$ is $\f 1e \theta_{\a_i\i\a_j}(q)$, which we
recall depends only on the ideal class of the fractional ideal
$\a_i\i\a_j$.  If $\chi$ is any character of
$\operatorname{Cl}(\O)$, viewed as the $h$-vector whose $i$-th
entry is $\chi(\a_i)$, then the $i$-th entry of the vector
$A_q\chi$ may be evaluated through matrix multiplication as
\begin{equation}  \label{aqchi}
(A_q\chi)(\a_i) \ = \  \f 1 e \sum_{\a_j \in\operatorname{Cl}(\O)}
\, \theta_{\a_i\i\a_j}(q)\,\chi(\a_j)  \ = \  \f 1 e \(
\sum_{\a_j\in \operatorname{ Cl}(\O)}\, \chi(\a_j)\,
\theta_{\a_j}(q) \) \chi(\a_i)\,,
\end{equation}
where in the last equality we have reindexed $\a_j\mapsto
\a_i\,\a_j$ using the group structure of $\operatorname{Cl}(\O)$.
Therefore $\chi$ is in fact an eigenvector of the matrix $e A_q$,
with eigenvalue equal to the sum of $\theta$-functions enclosed in
parentheses, known as a {\em Hecke $\theta$-function} (see
\cite[\S12]{iwaniec-blue}).  These, which we shall denote
$\theta_\chi(q)$, form a more natural basis of modular forms than
the ideal class $\theta$-functions $\theta_{\a_j}$ because they
are in fact  Hecke eigenforms. Using (\ref{genfunc}), the
$L$-functions of these Hecke characters can be written as
\begin{equation}  \label{heckelfunc}
\gathered
  L(s,\chi) \  \ =  \ \   L(s,\theta_\chi) \ \ = \ \
\sum_{\text{integral ideals } \mathfrak{a} \subset K}
\chi(\mathfrak{a})\,(N\mathfrak{a})^{-s} \ \ = \ \
\sum_{n=1}^\infty \, a_n(\chi)\, n^{-s}\,,
\\
\text{where}\qquad a_n(\chi) \ \ = \ \
\sum_{\overset{\scriptstyle{\text{integral ideals } \mathfrak{a}
\subset K}}{N\mathfrak{a} = n}} \chi({\mathfrak{a}} )
\qquad\qquad\qquad \qquad\qquad \ \
\endgathered
\end{equation}
is in fact simply the eigenvalue of $\,e\,M(n)\,$ for the
eigenvector formed from the character $\chi$ as above, which can
be seen by isolating the coefficient of $q^n$ in the sum on the
righthand side of~\eqref{aqchi}.

\subsection{Eigenvalue Separation under the Generalized Riemann Hypothesis}

Our isogeny graph is a superposition of the previous graphs
$M(n)$, where $n$ is a prime bounded by a parameter $m$ (which we
recall is $(\log q)^{2+\d}$ for some fixed $\d>0$). This
corresponds to a graph on the elliptic curves represented by ideal
classes in an order $\O$ of $K=\Q(\sqrt{d})$, whose edges
represent isogenies of prime degree $\le m$.   The graphs with
adjacency matrices $\{M(p) \mid p\le m\}$ above share common
eigenfunctions (the characters $\chi$ of $\operatorname{Cl}(\O)$),
and so their eigenvalues are
\begin{equation}\label{isoeig}
    \l_\chi \ \ = \ \ \f 1 e \,  \sum_{p\,\le\,m}\,a_p(\chi)\,
     \ \ = \ \ \f 1 e \,
    \sum_{p\,\le\,m}\, \sum_{\stackrel{\scriptstyle{\text{integral ideals }
     {\mathfrak a} \, \subset \, K}}{N{\mathfrak a} \, = \, p}}
    \chi({\mathfrak a})\,.
\end{equation}
When $\chi$ is the trivial character, $\l_{\text{triv}}$ equals
the degree of the regular graph $\Graph$.  Since roughly half of
rational primes $p$ split in $K$, and those which do split into
two ideals of norm $p$, $\l_{\text{triv}}$ is roughly
$\f{\pi(m)}{e} \sim \f{m}{e\log m}$ by the prime number theorem.
This eigenvalue is always the largest in absolute value, as can be
deduced from (\ref{isoeig}), because $|\chi({\mathfrak a})|$
always equals 1 when $\chi$ is the trivial character. For the
polynomial mixing of the random walk in Theorem~\ref{newmainthm}
we will require a separation between the trivial and nontrivial
eigenvalues of size $1/\text{polylog}(q)$. This would be the case,
for example, if for each nontrivial character $\chi$ there merely
exists one ideal $\mathfrak a$ of prime norm $\le m$ with
$\operatorname{Re}\chi({\mathfrak a}) \le 1-
\f{1}{\text{polylog}(q)}$.  This is analogous to the problem of
finding a small prime nonresidue modulo, say, a large prime $Q$,
where one merely needs to find any cancellation at all in the
character sum $\sum_{p\le m} \jacobi{p}{Q}$.  However, the latter
requires a strong assumption from analytic number theory, such as
the Generalized Riemann Hypothesis (GRH).  In the next section we
will accordingly  derive such bounds for $\l_\chi$, under the
assumption of GRH.  As a consequence of the more general
\lemref{filllater} we will show the following.

\begin{lem}\label{boundforchipsum}
Let  $D<0$ and let $\O$ be the quadratic order of discriminant
$D$. If $\chi$ is a nontrivial ideal class character of $\O$, then
the Generalized Riemann Hypothesis for $L(s,\chi)$ implies that
the sum (\ref{isoeig}) is bounded by $O(m^{1/2}\log |mD|)$ with an
absolute implied constant.
\end{lem}

\begin{proof}[of Theorem~\ref{newmainthm}]
There are
only finitely many levels for $q$ less than any given bound, so it
suffices to prove the theorem for $q$ large and $p(x)=x^{2+\d}$,
where $\d>0$ is fixed. The eigenvalues of the adjacency matrix for
a given level are given by (\ref{isoeig}). Recall that $|D|\le 4
q$ and $\l_{\hbox{triv}}\sim\f{m}{e\log m}$. With our choice of
$m=p(\log q)$, the bound for the nontrivial eigenvalues in
\lemref{boundforchipsum} is $\l_\chi=O(\l_{\text{triv}}^\b)$ for
any $\b>\f{1}{2}+\f{1}{\d+2}$.  That means indeed our isogeny
graphs are expanders for $q$ large; the random walk assertion
follows from this bound and Proposition~\ref{expander-prop}.
\end{proof}

\begin{proof}[of Corollary~\ref{newmaincor}]
The Theorem shows that a random walk from any fixed curve $E$
probabilistically reaches the proportion where the algorithm $A$
succeeds, in at most $\operatorname{polylog}(q)$ steps. Since each
step is a low degree isogeny, their composition can be computed in
$\operatorname{polylog}(q)$ steps.  Even though the degree of this
isogeny might be large, the degrees of each step are  small. This
provides the random polynomial time reduction of \textsc{dlog}
along successive curves in the random walk, and hence from $E$ to
a curve for which the algorithm $A$ succeeds.
\end{proof}

\section{The Prime Number Theorem for Modular Form $L$-functions}
\label{pntsec}

In this section we prove \lemref{boundforchipsum}, assuming the
Generalized Riemann Hypothesis (GRH) for the $L$-functions
(\ref{heckelfunc}).  Our argument is more general, and in fact
gives estimates for sums of the form $\sum_{p\le m} a_p$, where
$a_p$ are the prime coefficients of any $L$-function.  This can be
thought of as an analog of the Prime Number Theorem because for
the simplest $L$-function, $\zeta(s)$, $a_p=1$ and this sum is in
fact exactly $\pi(m)$.  As a compromise between readability and
generality, we will restrict the presentation here to the case of
modular form $L$-functions (including (\ref{heckelfunc})).
Background references for this section include
\cite{iwaniec,iwaniec-blue,murty}; for information about more
general $L$-functions see also \cite{gelbart-miller,rudsar}.

We shall now consider a classical holomorphic modular form $f$,
with Fourier expansion  $f(z)=
    \sum_{n\,=\,0}^\infty\,c_n\,e^{2\pi i n z}$.
We will assume that $f$ is a Hecke eigenform, since this condition
is met in the situation of  \lemref{boundforchipsum} (see the
comments between (\ref{aqchi}) and (\ref{heckelfunc})). It is
natural to study the renormalized coefficients $a_n =
n^{-(k-1)/2}c_n$, where $k\ge 1$ is the weight of $f$ (in
\secref{thetsubsec} $k=1$, so $a_n=c_n$). The $L$-function
 of such a
 modular form can be written as the Dirichlet series
  $ L(s,f) = \sum_{n=1}^\infty  a_n   n^{-s}   =
\prod_p\, (1-\a_{p} p^{-s})^{-1} (1-\b_{p} p^{-s})^{-1}$, the last
equality using the fact that  $f$ is a Hecke eigenform. The
$L$-function $L(s,f)$ is entire when $f$ is a cusp form
(e.g.~$a_0=0$).  The Ramanujan conjecture (in this case a theorem
of \cite{deligne} and \cite{delser}) asserts that $ \, |\a_p|,
|\b_p| \, \le \, 1$.

\lemref{boundforchipsum} is  concerned with estimates for the sums
\begin{equation}\label{smfdef}
    S(m,f)   \ \ :=
 \  \  \sum_{p \, \le \,  m}  \ a_p\, .
\end{equation}
  As with the prime number
theorem, it is more convenient to instead analyze the weighted sum
\begin{equation}\label{weightedsum}
\psi(m,f)   \ \ :=  \ \   { \sum_{p^k}}  \ b_{p^k}  \, \log p
\end{equation}
over prime powers, where the coefficients $b_n$ are those
appearing in the Dirichlet series for $-\f{L'}{L}(s)$:
$$ -\,{  \f{L'}{L}}(s)  \ \  = \ \ { \sum_{n\,=\,1}^\infty}\
  b_n  \, \L(n) \, n^{-s}
\ \  =  \ \  { \sum_{p,\, k} } \ b_{p^k} \,  \log(p) \,
p^{\,-k\,s}\,,$$ i.e., $b_{p^k}  \ =  \   \a_p^k  \  +  \ \b_p^k$.

\begin{lem}  For a holomorphic modular form $f$ one has
$$ \ \psi(m,f)  \ = \   \sum_{p \, \le \,  m}  a_p  \log p \ +
 \ O(m^{1/2}).$$
\end{lem}

\begin{proof}  The error term represents the contribution of proper prime
powers.  Since $|b_{p^k}|\le 2$, it is bounded by twice
\begin{equation}\label{bddabove}
\sum_{{\stackrel{\scriptstyle{p^k \, \le \,  m}}{ k \, \ge \, 2}}}
\log p \ \ = \ \  \sum_{\stackrel{\scriptstyle{p \, \le \,
m^{1/2}}}{ 2  \, \le \,  k \, \le \,  \f{\log m}{\log p}}} \log p
\ \  \le \ \  \sum_{p \, \le \,  m^{1/2}}  \, \log p  \, \f{\log
m}{\log p} \ \ \le \ \ \pi(m^{1/2})\log m \, ,
\end{equation}
which is $O(m^{1/2})$ by the Prime Number Theorem.
\end{proof}

\begin{lem}(Iwaniec~\cite[p. 114]{iwaniec})
 Assume that $f$ is a holomorphic modular cusp form of level\footnote{Actually in \cite{iwaniec} $N$ equals the conductor of
the $L$-function, which in general may be smaller than the level.
The lemma is of course nevertheless valid.} $N$ and that  $L(s,f)$
satisfies GRH. Then $\psi(m,f)  =  O(m^{1/2} \log(m)
\log(m N))$.
\end{lem}

We deduce that $
  S'(m,f) :=   {\sum_{p\le m}} a_p  \log p =    O(m^{1/2} \log(m)
\log(mN))$.   Finally we shall estimate the sums $S(m,f)$ from
(\ref{smfdef}) by removing the $\log(m)$ using a standard partial
summation argument.

\begin{lem}
\label{filllater} Suppose that $f$ is a holomorphic modular cusp
form of level $N$ and $L(s,f)$ satisfies GRH.  Then
$S(m,f)=O(m^{1/2}\log (mN))$.
\end{lem}
\begin{proof}
 First  define $\tilde{a}_p$ to be $a_p$,
if $p$ is prime, and 0 otherwise.  Then
$$
\sum_{p \, \le \,  m} a_p  \ \  = \ \  \sum_{p \, \le \,  m}\,
[\tilde{a}_p\log p]\,\f{1}{\log p} \ \   = \ \  \sum_{n\, \le \,
m}\,[\tilde{a}_n\log n]\,\f{1}{\log n}\,.
$$
By partial summation over $2\le n \le m$, we then find
\begin{eqnarray*}
\sum_{p\,\le \, m}\,a_p & =  & \sum_{n\,<\,m} \,S'(n,f)
\left( \f{1}{\log (n)} \, - \, \f{1}{\log(n+1)}\right) \ + \ \f{S'(m,f)}{\log m} \\
& \ll & \sum_{n\,<\, m}\, \(n^{1/2}\log(n)\log(nN)\)\,
 \left| \f{d}{dn}\left( ( \log n  )^{-1}\right)\right| \ + \ m^{1/2}\,\log(mN)\\
& \ll &  \sum_{n \, < \,  m} n^{1/2} \log(n) \log(nN) \f{1}{n(\log
n)^2} \ + \ m^{1/2}\,\log(mN)\, ,
\end{eqnarray*}
so in fact $S(m,f) = \sum_{p \le  m} \, a_p  = O( m^{1/2} \log
(mN))$.
\end{proof}

All the implied constants in these 3 lemmas are absolute.
 Some useful estimates for them may be found in \cite{bach}.

\subsection{Subexponential Reductions via Lindel\"of Hypothesis}
\label{lindelof} In the previous lemma we have assumed GRH.  It
seems very difficult to get a corresponding unconditional bound
for $S(m,f)$. However, a slightly weaker statement can be proven
by assuming only the Lindel\"of hypothesis (which is a consequence
of GRH). Namely, one has that $\sum_{n\le m }\,a_n =
O_\e(m^{1/2+\e}N^\e)$, for any $\e>0$
(\cite[(5.61)]{iwaniec-blue}).  The fact that this last sum is
over all $n\le m$, not just primes, is not of crucial importance
for our application.  However, the significant difference here is
that the dependence on $N$ is not polynomial in $\log N$, but
merely subexponential.  This observation can be used to weaken the
hypothesis in Theorem~\ref{newmainthm} and
Corollary~\ref{newmaincor} from GRH to the Lindel\"of hypothesis,
at the expense of replacing ``polynomial'' by ``subexponential.''

\section{Distribution of $c_\pi$} \label{distr}

Theorem~\ref{newmainthm} and Corollary~\ref{newmaincor} are statements
about individual levels.  As we mentioned in \secref{summary}, our
random reducibility result extends between two levels as long as the
levels satisfy the requirement that their conductors differ by
polynomially smooth amounts. In this section we explore this extension
in more detail, and explain why the above requirement is typically
satisfied.

It was mentioned after Theorem~\ref{direction-classification} that
the largest possible conductor difference is $c_\pi$, which is the
largest square factor of $d_\pi= \Trace(E)^2 - 4q$.  In principle
this factor could be as large as $2\sqrt{q}$, though statistically
speaking most integers (a proportion of $\f{6}{\pi^2}\approx.61$)
are square-free, explaining why $c_\pi$ is very often 1 or at
least fairly small \cite{tenenbaum}.  This means, for example,
that most randomly selected elliptic curves have an isogeny class
consisting of only one level.

When an isogeny class consists of multiple levels, we need to be
able to construct vertical isogenies between levels in order to
conclude that \textsc{dlog} instances between the levels are
randomly reducible to each other. The fastest known algorithm for
constructing vertical isogenies between two levels, due to
Kohel~\cite{kohel}, has runtime $O(\ell^4)$, where $\ell$ is the
largest prime dividing the conductor of one of the levels, but not
the other.  Any two levels which can be efficiently bridged via
Kohel's algorithm can be considered as one unit for the purposes
of random reducibility. Accordingly, polynomial time random
reducibility holds within an isogeny class if $c_\pi$ for that
isogeny class is polynomially smooth.

With this in mind, we will now determine a heuristic estimate for
the expected size of the largest prime factor $P(c_\pi)$ of
$c_\pi$, i.e., the largest prime which divides $d_\pi$ to order at
least $2$.  The trace $t=\Trace(E)$, when sampled over random
elliptic curves, is thought to have a fairly uniform distribution
over most of the Hasse interval.  This serves to predict the
useful heuristic that $-d_\pi=4q-t^2$ is typically of size $q$
(see for example \cite{lenstra,serre}).  Assuming that, the
probability that $P(c_\pi)$ exceeds $\b$ can be loosely estimated
as $O(1/\b)$. This is because roughly a fraction of $\rho =
\prod_{p\,>\,\b}^{\sqrt{q}}\(1\,-\,p^{-2}\)$ integers of size $q$
have no repeated prime factor $p>\b$.  It is easy to see that
$\log(\rho)=O(\sum_{n>\b}n^{-2})=O(1/\b)$, so that
$1-\rho=O(1/\b)$ as suggested.

It follows that a randomly selected elliptic curve is extremely likely
to have a small enough value of $P(c_\pi)$ to allow for random
reducibility throughout its entire isogeny class.  This explains why
in Figure~\ref{standardscurves} all of the randomly generated curves
have $P(c_\pi) = 1$, except for one curve which has $P(c_\pi) = 3$.


Finally, let us  consider the situation where a {\em non-random}
curve is deliberately selected so as to have a large value of
$c_\pi$.  Currently the only known methods for constructing such
curves is to use complex multiplication methods~\cite[Ch.
VIII]{blake} to construct curves with a predetermined number of
points chosen to ensure that $c_\pi$ is almost as large as
$\sqrt{d_\pi}$. Some convenient examples of such curves are the
Koblitz curves listed in the NIST FIPS 186-2
document~\cite{standardsdocumentsnist}, which we have also
tabulated in Figure~\ref{standardscurves}.  Since these curves all have
complex multiplication by the field $K=\Q(\sqrt{-7})$, the
discriminants of these curves are of the form $d_\pi=-7 c_\pi^2$.
If we assume that $c_\pi$ behaves as a random integer of size
$\sqrt{d_\pi}$, which is roughly $\sqrt{q}$, then the distribution
of $P(c_\pi)$ is governed by the usual smoothness bounds for large
integers \cite{tenenbaum}, and hence is typically too large to
permit efficient application of Kohel's algorithm for navigating
between levels. Thus we cannot prove random reducibility from a
theoretical standpoint for all of the elliptic curves within the
isogeny class $S_{N,q}$ of such a specially constructed curve.
However, in practice only a small subset of the elliptic curves in
$S_{N,q}$ are efficiently constructible using the complex
multiplication method (or any other presently known method), and
this subset coincides exactly with the subcollection of levels in
$S_{N,q}$ which are accessible from the top level (where
$\End(E)=\mathcal O_K$) using Kohel's algorithm. Pending future
developments, it therefore remains true that all of the special
curves that we can construct within an isogeny class
have equivalent \textsc{dlog} problems in the random reducible
sense.

\vspace{.2cm}

 {\bf Acknowledgments:} It is a pleasure to thank
William Aiello, Michael Ben-Or, Dan Boneh, Brian Conrad, Adolf
Hildebrand, Henryk Iwaniec, Dimitar Jetchev, Neal Koblitz,
Alexander Lubotzky, Peter Sarnak, Adi Shamir, and Yacov Yacobi for
their discussions and helpful comments. We are also indebted to
Peter Montgomery for his factoring assistance in producing
Figure~\ref{standardscurves}.

\appendix

\section{Supersingular Case}

\label{ss-sec}

In this appendix we discuss the isogeny graphs for supersingular
elliptic curves and prove Theorem~\ref{newmainthm} in this
setting. The isogeny graphs were first considered by
Mestre~\cite{mestre}, and were shown by Pizer~\cite{pi90,piz98} to
have the Ramanujan property.  Curiously, the actual graphs were
first described by Ihara \cite{ihara} in 1965, but not noticed to
be examples of expander graphs until much later.  We have decided
to give an account here for completeness, mainly following Pizer's
arguments. The isogeny graphs we will present here differ
from those in the ordinary case in that they are \emph{directed}.
This will cause no serious practical consequences, because one can
arrange that only a bounded number of edges in these graphs will
be unaccompanied by a reverse edge. Also, the implication about
rapid mixing used for Theorem~\ref{newmainthm} carries over as
well in the directed setting with almost no modification. It is
instructive to compare the proofs for the ordinary and
supersingular cases, in order to see how GRH plays a role
analogous to the Ramanujan conjectures.

Every $\bar{\F}_q$-isomorphism class of supersingular elliptic curves
in characteristic $p$ is defined over either $\F_p$ or
$\F_{p^2}$~\cite{silverman}, so it suffices to fix $\F_q=\F_{p^2}$ as
the field of definition for this discussion. Thus, in contrast to
ordinary curves, there is a finite bound $g$ on the number of
isomorphism classes that can belong to any given isogeny class (this
bound is in fact the genus of the modular curve $X_0(p)$, which is
roughly $\f{p+1}{12}$). It turns out that all isomorphism classes of
supersingular curves defined over $\F_{p^2}$ belong to the same
isogeny class~\cite{mestre}. Because the number of supersingular
curves up to isomorphism is so much smaller than the number of
ordinary curves up to isomorphism, correspondingly fewer of the edges
need to be included in order to form a Ramanujan graph.  For a fixed
prime value of $\ell \neq p$, we define the vertices of the
supersingular isogeny graph $\Graph$ to consist of these $g$
isomorphism classes, with directed edges indexed by equivalence
classes of degree-$\ell$ isogenies as defined below. In fact, we will
prove that $\Graph$ is a directed $k=\ell+1$-regular graph
satisfying the Ramanujan bound of $|\lambda| \leq 2
\sqrt{\ell}=2\sqrt{k-1}$ for the nontrivial eigenvalues of its
adjacency matrix.  The degree $\ell$ in particular may be taken to be
as small as 2 or 3.

For the definition of the equivalence classes of isogenies --- as well
as later for the proofs --- we now need to recall the structure of the
endomorphism rings of supersingular elliptic curves.  In contrast to
the ordinary setting (\secref{prelim}), the endomorphism ring
$\End(E)$ is a maximal order in the quaternion algebra
$R=\Q_{p,\infty}$ ramified at $p$ and $\infty$.  Moreover, isomorphism
classes of supersingular curves $E_i$ isogenous to $E$ are in 1-1
correspondence with the left ideal classes $I_i := \Hom(E_i, E)$ of
$R$. As in Section~\ref{graph-specification}, call two isogenies
$\phi_1, \phi_2\colon E_i \to E_j$ equivalent if there exists an
automorphism $\alpha$ of $E_j$ such that $\phi_2 = \alpha
\phi_1$. Under this relation, the set of equivalence classes of
isogenies from $E_i$ to $E_j$ is equal to $I_j^{-1} I_i$ modulo the
units of $I_j$. This correspondence is degree preserving, in the sense
that the degree of an isogeny equals the reduced norm of the
corresponding element in $I_j^{-1} I_i$, normalized by the norm of
$I_j^{-1} I_i$ itself.  This is the notion of equivalence class of
isogenies referred to in the definition of $\Graph$ in the previous
paragraph.  Thus, for any integer $n$, the generating function for the
number $M_{ij}(n)$ of equivalence classes of degree $n$ isogenies from
$E_i$ to $E_j$ (i.e.,~the number of edges between vertices representing
elliptic curves $E_i$ and $E_j$) is given by
\begin{equation}\label{ss-q-expansion}
\sum_{n=0}^\infty \, M_{ij}(n) \, q^n \ \  :=  \ \ \frac{1}{e_j}
\sum_{\alpha \, \in  \, I_j^{-1} I_i} q^{\,N(\alpha)/N(I_j^{-1}
I_i)}\, ,
\end{equation}
where $e_j$ is the number of units in $I_j$ (equivalently, the
number of automorphisms of $E_j$).  One knows that $e_j \le 6$,
and in fact $e_j=2$ except for at most two values of $j$ -- see
the further remarks at the end of this appendix. Proofs for the
statements in this paragraph can be found in~\cite{gr1,piz98}.

The $\theta$-series on the righthand side of
\eqref{ss-q-expansion} is a  weight 2 modular form for the
congruence subgroup $\Gamma_0(p)$, and the matrices
$$
B(n) :=
\begin{pmatrix}
M_{11}(n) & \cdots & M_{1g}(n) \\
\vdots & \ddots & \vdots \\
M_{g1}(n) & \cdots & M_{gg}(n)
\end{pmatrix}
$$ (called Brandt matrices) are simultaneously both the $n$-th Fourier coefficients
of various modular forms, as well the adjacency matrices for the
graph $\Graph$.  A fundamental property of the Brandt matrices
$B(n)$ is that they represent  the action of the $n^\mathrm{th}$
Hecke operator $T(n)$ on a certain basis of modular forms of
weight 2 for $\Gamma_0(p)$ (see \cite{pi90}). Thus the eigenvalues
of $B(n)$ are given by the $n^\mathrm{th}$ coefficients of the
weight-2 Hecke eigenforms for $\Gamma_0(p)$. These eigenforms
include a single Eisenstein series, with the rest being cusp
forms.  Now we suppose that $n=\ell$ is prime (mainly in order to
simplify the following statements).  The $n^\mathrm{th}$ Hecke
eigenvalue of the Eisenstein series is $n+1$, while those of the
cusp forms are bounded in absolute value by $2\sqrt{n}$ according
to the Ramanujan conjectures (in this case a theorem of Eichler
\cite{eichler} and Igusa \cite{igusa}). Thus the adjacency matrix
of $\Graph$ has trivial eigenvalue equal to $\ell+1$ (the degree
$k$), and its nontrivial eigenvalues indeed satisfy the Ramanujan
bound $|\l| \le 2\sqrt{k-1}$.

Finally, we conclude with some comments about the potential
asymmetry of the matrix $B(n)$.  This is due to the asymmetry in
the definition of equivalence classes of isogenies. Indeed, if
$\Aut(E_1)$ and $\Aut(E_2)$ are different, then two isogenies $E_1 \to
E_2$ can sometimes be equivalent even when their dual isogenies are
not equivalent. This problem arises only if one of the curves $E_i$
has complex multiplication by either $\sqrt{-1}$ or $e^{2 \pi
i/3}$, since otherwise the only possible automorphisms of $E_i$
are the scalar multiplication maps $\pm 1$~\cite[\S
III.10]{silverman}. In the  supersingular setting, one can avoid
curves with such unusually rich automorphism groups by choosing a
characteristic $p$ which splits in both $\Z[\sqrt{-1}]$ and
$\Z[e^{2 \pi i/3}]$, i.e.,~$p \equiv 1 \bmod 12$ (see \cite[Prop.
4.6]{pi90}). In the case of ordinary curves, however, the
quadratic orders $\Z[\sqrt{-1}]$ and $\Z[e^{2 \pi i/3}]$ both have
class number $1$, which then renders the issue moot because  the
isogeny graphs corresponding to these levels each have only one
node.

\end{document}